# Some Notes on Geometric Interpretation of Holding-Free Solution for Urban Intersection Model

Yen-Hsiang TanChen [a,1]

Feb, 2017

**Abstract**

Conventionally, methods to solve macroscopic node model are discussed in algorithm or in algebraic point of view. In this paper, the geometric interpretation is discussed, focusing on an example in Flötteröd and Rohde (2011). Via observing these diagrams, it is easy to see that deriving holding-free solution (HFS) is not just an ordinary optimisation problem, but rather a Pareto frontier-finding problem.

## 1. Introduction

For dynamic network loading (DNL) problems or other kinds of network in transportation system, there are two types of macroscopic traffic model in general: the link models and the node models. The former capture the bulk behaviour of road traffic in homogeneous sections; the latter describes dynamics at bottlenecks or junctions, which are point-like positions that connect multiple links and has certain properties that affect the evolution of states of incoming links and outgoing links.

The node model in general (in addition to simple merge and diverge on motorways or expressways) had not been explored thoroughly until recent discovery of many properties: invariance principle, uniqueness, Pareto optimality...etc. In Lebacque and Khoshyaran (2005), the two invariance principles are argued and only some of the models fulfil this property. The known models obeying these properties at general nodes are: (Flötteröd and Rohde, 2011, Tampere et al., 2011, Gibb, 2011, Courout et al., 2012, Smits et al., 2015, and Jabari, 2016). The uniqueness is discussed in (Courout et al., 2012). Pareto optimality is discussed in (Smits et al., 2015). The Pareto optimality must be satisfied if the model takes into account the individual flow maximisation, which states that no drivers would leave the capacity unused if they have chance to use it. Among the contributions in these articles, we mainly focus on the relation between Pareto optimality and holding-free solution (mentioned in the next paragraph) in the classic example in (Flötteröd and Rohde, 2011).

[1] Contact: b94901007@ntu.edu.tw (This email is permanent and does not reflect his current affiliation)
a) Greenwave Octopus Transport Solutions

Recently, condition of holding-free solution (HFS) is discussed in (Jabari, 2016) Jabari (2016) proves that the Incremental Node Transfer Model (INM) proposed in Flötteröd and Rohde (2011), which is also equivalent to other capacity proportional model like Tampere et al. (2011) (Smit et al., 2015), is indeed a HFS. It further states that in terms of flow maximisation, the solution of INM is suboptimal. However, does the optimality in terms of total flow necessary? If not, what are the steps in INM do? Jabari (2016) also show that the flow maximising and a greedy algorithm is HFS. What are common properties of them, then?

We investigate about these questions from observing properties of geometry in space form by flows. The rest of the article is organised as follows, we first review on how to draw diagrams by hand in Section 2. The diagrams of Flötteröd-Rohde Example and what different algorithms do in flow space are shown in Section 3. Inspired by diagrams, one can infer the diagrams are actually exhibit Pareto optimal, shown in Section 4 with another way of formulation of HFS. In Section 5, some topics discussed and summary is made in Section 6.

## 2. Plane (Hyperplane) and Intercept

Demand and supply constraints form inequality in macroscopic node modelling. In general, constrains from demand and supply are hyperplanes in space. However, in the Flötteröd- Rohde Example, there are only 3 incoming links, thus drawing on a 3-dimensional space is possible.

Rearrange into *intercept form* one can directly find intercepts at axes.

$$\frac{x_1}{a_1}+\frac{x_2}{a_2}+\frac{x_3}{a_3}=1 \qquad (1)$$

where $a_i$ is the $x_i$-axis intercept. Since all the demand and supply constraints can be rearranged into this form, one can easily draw the diagram by hand.

The feasible region is in general polyhedral; However, some models having *internal supply constraints* (for example, models in Flötteröd and Rohde, 2011, and Tampere et al., 2011) introduce complicates the node. If with the additional *internal supply constraints,* (It correspond to the term "potential capacity" formula in Highway Capacity Manual 2011) the feasible region would not be polyhedral. Thus, in this paper, we focus only on node problem without *internal supply constraints.*

For INM algorithm, if in any step the supply constraint is binding (active), then merging weights must be checked between competing incoming flows. The larger the merging weight (defined originally in Flötteröd and Rohde (2011)) $\boldsymbol{\alpha}=[\alpha_1, \alpha_2, \alpha_3]^T$, the much more flow of an incoming link can utilise the residual supply (within a step of incremental). In geometry, it is actually a line in space, with the following formula:

$$\mathbf{L}: \quad (q_1^{(k-1)} + p \cdot z_1^{(k)} \cdot \alpha_1, \qquad \ldots, q_i^{(k-1)} + p \cdot z_i^{(k)} \cdot \alpha_i \ldots) \tag{2}$$

where, k is the iteration step,

p: $p \in (-\infty, \infty)$ parameter describing this line

$z_i^{(k)}$: binary indicator indicating whether $\{f(i,o)>0\} \wedge \{ i \text{ in } \mathcal{D}_{in}\}$

$\mathcal{D}_{in}$: the set of incoming links that are still within consideration of algorithm of INM algorithm in (Flötteröd and Rohde, 2011)

The above simple geometry fundamentals allow us to draw algorithm steps in flow space by hand easily.

## 3. The Flötteröd- Rohde Example

### 3.1 Problem Statement

An example appeared firstly in Flötteröd and Rohde(2011) can be described as**:** A node consists of 3 incoming links (in $I$, from west, north, and south, respectively) and 3 outgoing links (in $\mathcal{O}$, to north, west, and south, respectively). The flow vector is defined as: $\mathbf{q} \equiv [q_1, q_2, q_3]^T \equiv [q_E, q_N, q_S]^T$ and they have demand constraints of, $\delta_1$=100 , $\delta_2$=600,and $\delta_3$=600, respectively. The supply constraints of outgoing links are: $\sigma_4 \equiv \sigma_N$=1400, $\sigma_5 \equiv \sigma_W$=400, and $\sigma_6 \equiv \sigma_S$=1400, respectively. Turing ratios from i to o are: $f(i=1, o=5)$ = 1.0, $f(2,5)$= $f(2,6)$ = 0.5, and $f(3,4)$= $f(3,5)$ = 0.5.

In the following sections, let us investigate its geometry properties first and then see what various algorithms are doing in their computation procedure.

### 3.2 Demand Constraint and Nature of Flow

The demand constraints and nature of flow makes result in a feasible box in the flow space, depicted in Fig 1. The demand constraints in Flötteröd- Rohde Example are:

$$\begin{cases} \mathbf{\Gamma_1}: & q_1 \leq \delta_1 = 100 \\ \mathbf{\Gamma_2}: & q_2 \leq \delta_2 = 600 \\ \mathbf{\Gamma_3}: & q_2 \leq \delta_3 = 600 \end{cases} \tag{3}$$

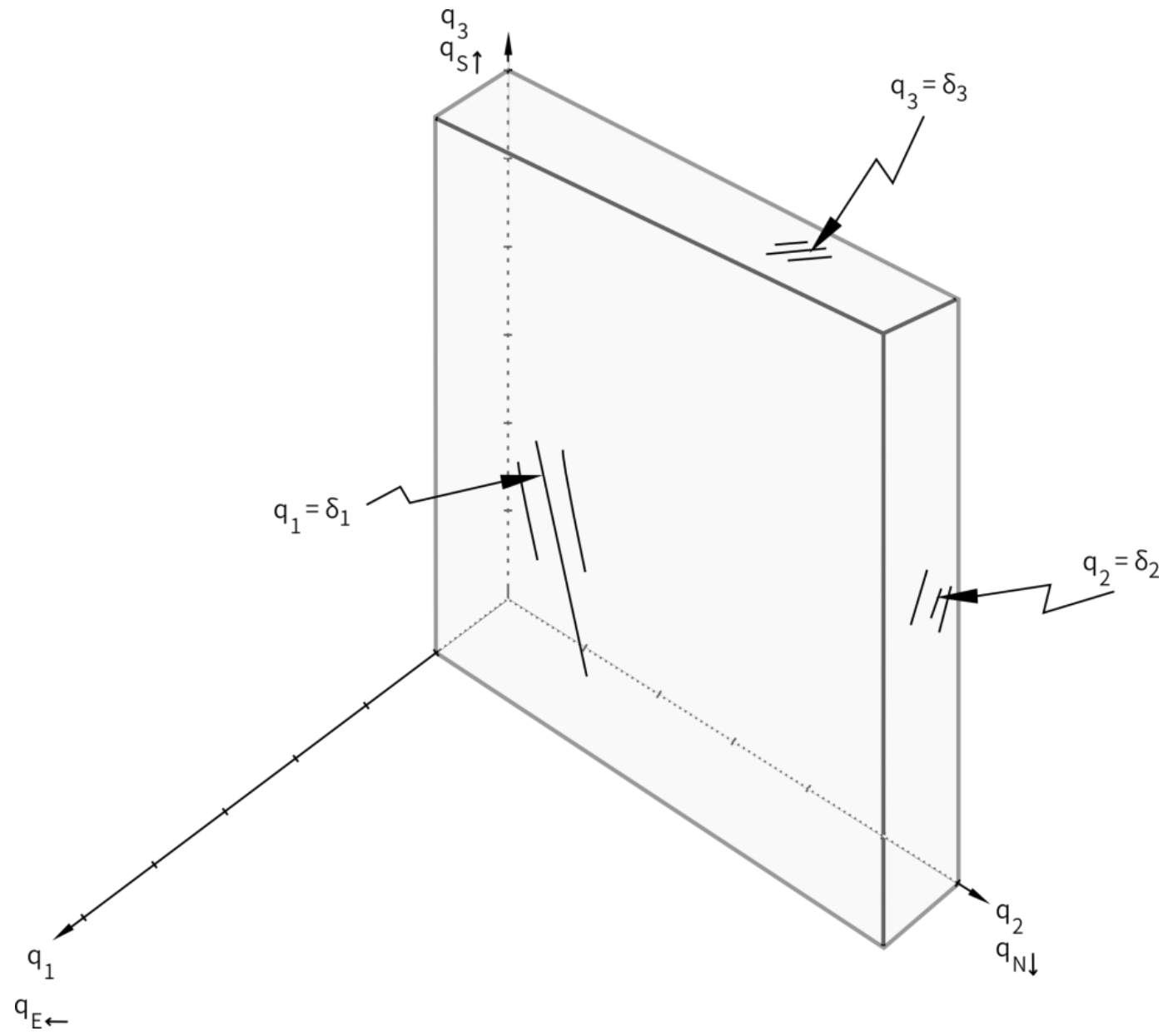

**Fig 1** Feasible box form by demand constraints and nature of flow.

The nature of flow (non-negativity, Tampere et al., 2011)forms the (trivial) inequality of:

$$\boldsymbol{\Gamma_0}: \begin{cases} q_1 \geq 0 \\ q_2 \geq 0 \\ q_3 \geq 0 \end{cases} \tag{4}$$

**3.2 Supply Constraints**

The supply constraints are:

$$\begin{gathered} 0.5q_3 \leq \sigma_4 = 1400 \\ q_1 + 0.5q_2 + 0.5q_3 \leq \sigma_5 = 400 \\ 0.5q_2 \leq \sigma_6 = 1400 \end{gathered} \tag{5}$$

If supply constraints are rearranged in intercept form, they appear like:

$$\begin{cases} \boldsymbol{\Gamma_4}: & q_3/2800 \leq 1 \\ \boldsymbol{\Gamma_5}: & \frac{q_1}{400} + \frac{q_2}{800} + \frac{q_3}{800} \leq 1 \\ \boldsymbol{\Gamma_6}: & q_2/2800 \leq 1 \end{cases} \tag{6}$$

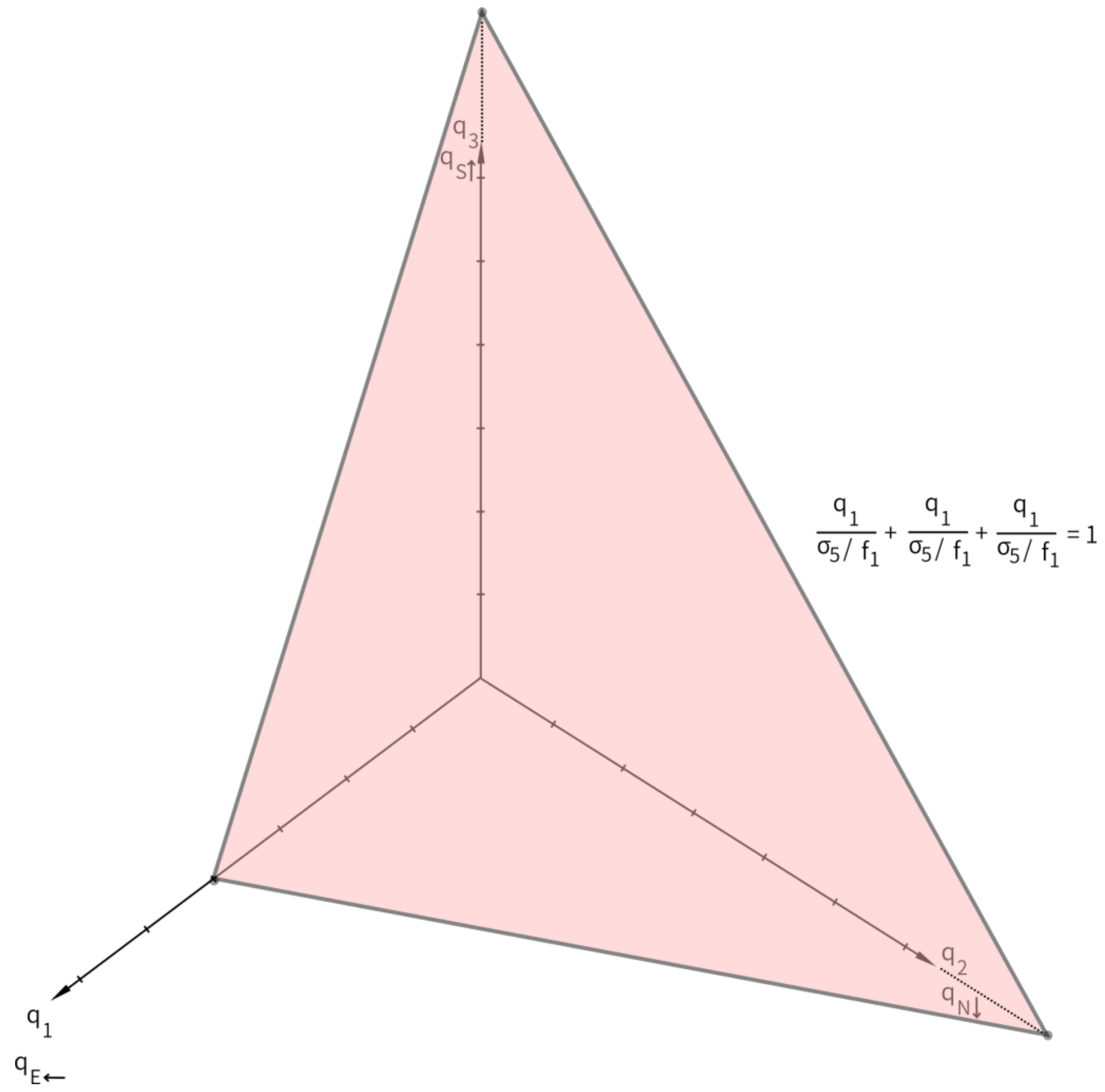

**Fig 2** Feasible region form by west outgoing arm (arm 5)

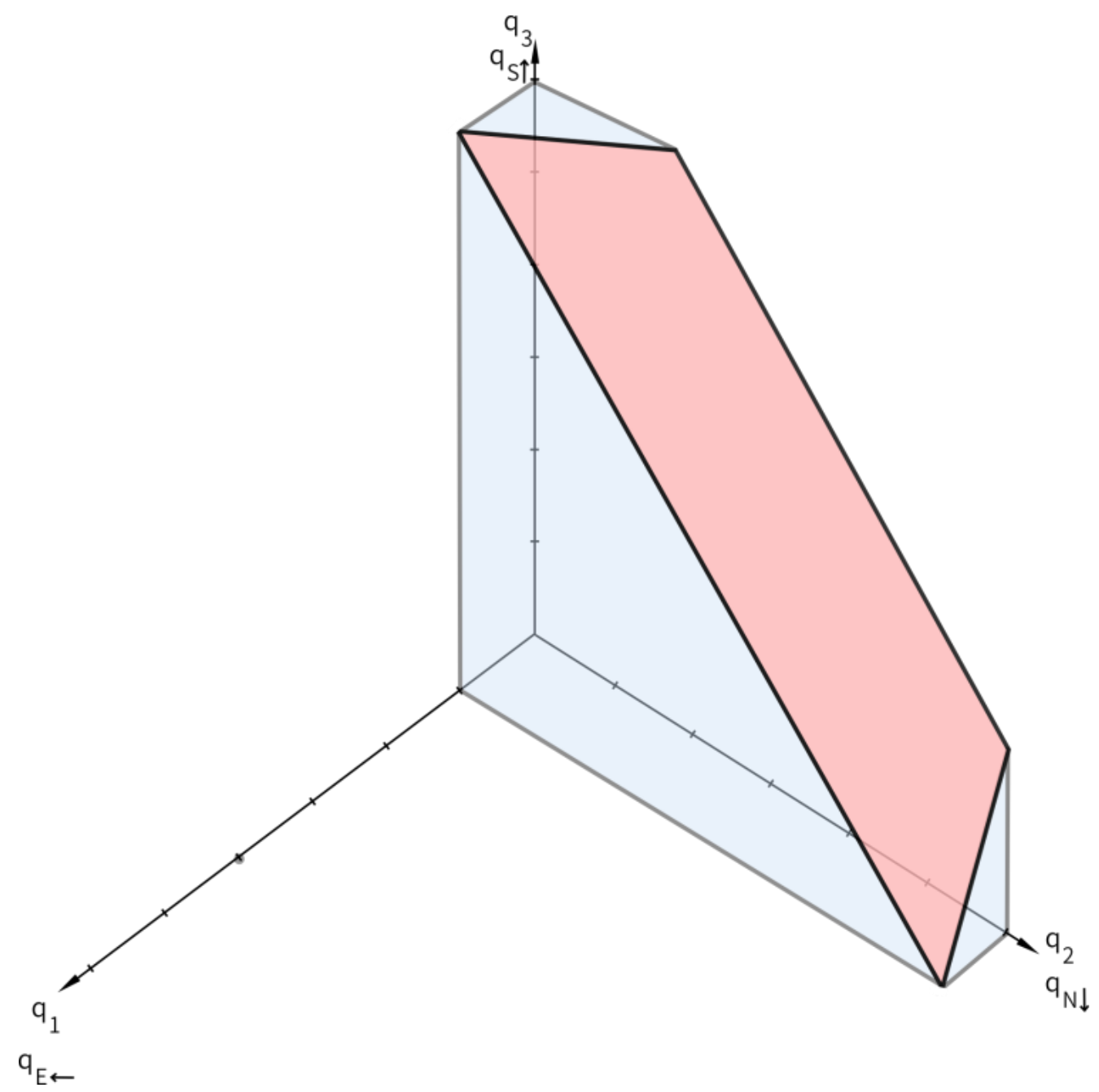

**Fig 3** The polyhedral feasible region.

### 3.3 Supply and Demand

In Fig 2, only the west outgoing arm supply constraint ($\boldsymbol{\Gamma}_5$) is depicted because other outgoing arms are inactive constraint due to their axis intercepts are larger than that form by demand constraints.

Considering the intersect of both supply and demand constraints (and trivially the positive flow nature constraint), the feasible region has a polyhedral shape as shown in Fig 3.

### 3.4 Incremental Node Transfer

The solution from INM (complete algorithm in Flötteröd and Rohde (2011)) is $\mathbf{q}_{\mathrm{INM}}^{*}$ =[16.7, 600, 166.7]$^{\mathrm{T}}$, which can be seen in Fig 4. With permutation of $\mathcal{P}_{\mathrm{INM}}$=<$\mathrm{i_N}$, $\mathrm{i_E}$, $\mathrm{i_s}$> and merging weight of $\boldsymbol{\alpha}$ =[0.1, 10, 1]$^{\mathrm{T}}$. The steps go like this in geometry point of view:

**Step 1:** Move to binding condition (when equality holds in that equality) of $\Gamma_2$. It is point A in Fig 4.

**Step 2a:** Start from [0, 600, 0]$^{\mathrm{T}}$, now let us try the extension to the parallel direction of axis $\mathrm{q}_2$ (direction of $\mathrm{i_E}$). We encounter a supply constraint (binding of $\Gamma_5$), thus the merging weight must be taking into consideration.

**Step 2b:** The line corresponds to merging weight is:

$$\mathbf{L}^{(\mathrm{k}=2)}: \quad (0 + \mathrm{p} \cdot 1 \cdot (0.1),\ 600 + \mathrm{p} \cdot 0 \cdot (10),\ 0 + \mathrm{p} \cdot 1 \cdot 1), \tag{7}$$

$$\mathrm{p} \in (-\infty, \infty)$$

It is the bold dashed line ($\overleftrightarrow{\mathrm{A}\,\mathbf{q}_{\mathrm{INM}}^{*}}$) shown in Fig 5. The algorithm terminates here, and the termination point is the intersection of the line and the (hyper) plane form by supply constraint:

$$\mathbf{q}_{\mathrm{INM}}^{*}: \begin{cases} \text{binding of } \Gamma_5 \\ \mathrm{L}^{(\mathrm{k}=2)} \end{cases} \tag{8}$$

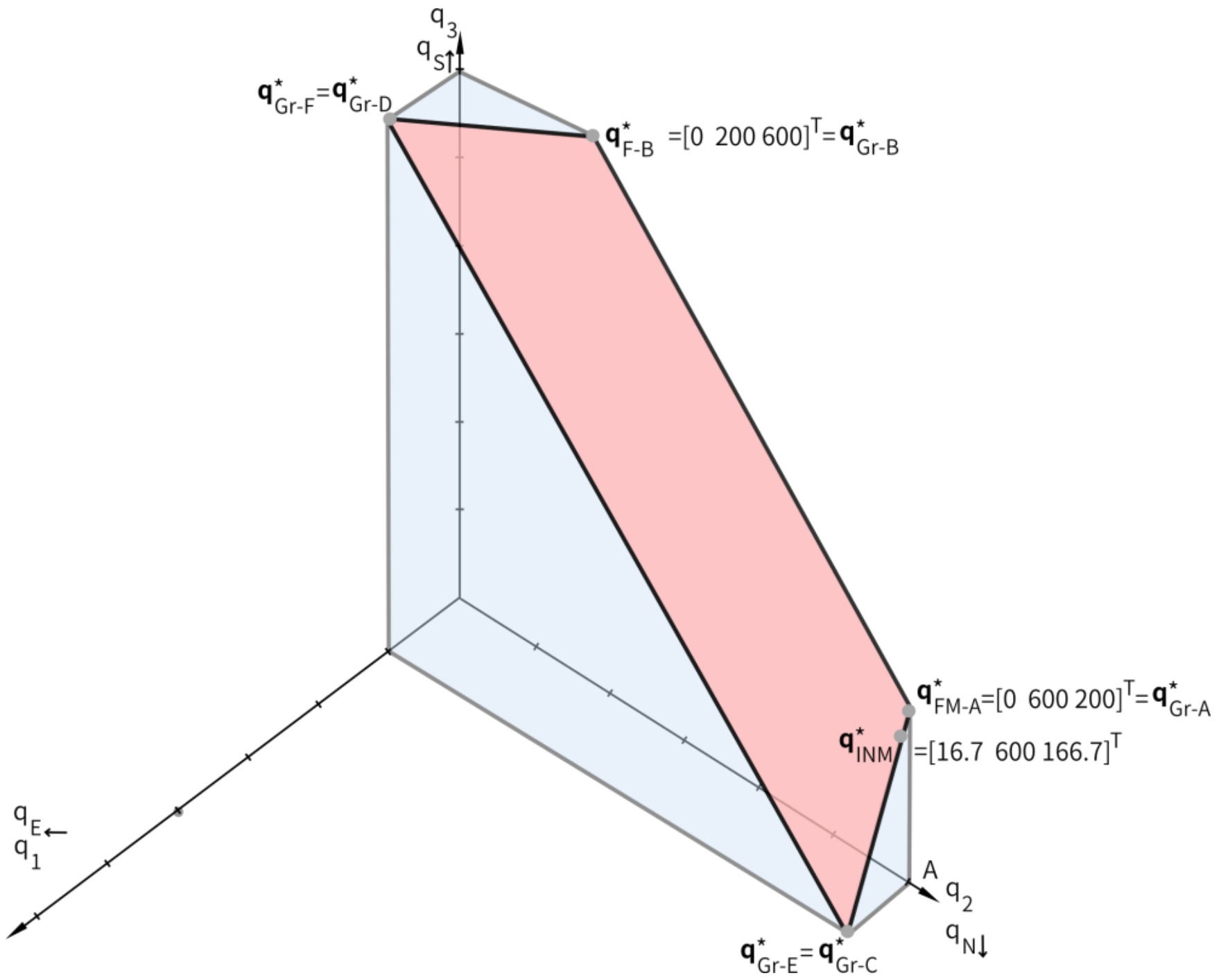

**Fig 4** The solutions from different algorithm, including INM, FM-x (from flow maximisation), and Gr-x (from greedy algorithm with various permutations).

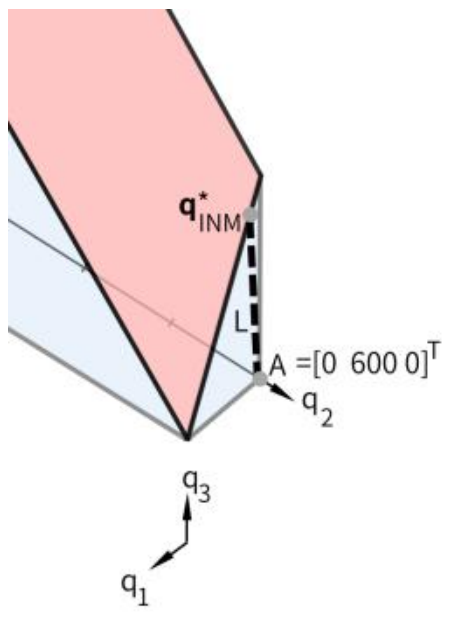

**Fig 5** Step 1 of INM steps on point A, and knowing that the (hyper)plane is binding on $q_3$ direction (and also the $q_1$ direction) to the next step, an auxiliary line L can be drawn to intersect with it, then the solution is found.

### 3.5 Flow Maximisation

The flow maximisation in Jabari (2016) can be arranged as:

$$\max \mu \tag{9}$$

s.t.

$$\begin{cases} \boldsymbol{\Lambda}: & \dfrac{\sum q_i}{\mu} = 1 \\ \boldsymbol{\Gamma_0}, \boldsymbol{\Gamma_1}, \boldsymbol{\Gamma_2} \dots \boldsymbol{\Gamma_6} & \end{cases} \tag{10}$$

In flow space, $\boldsymbol{\Lambda}$ is a (hyper) plane with normal direction of $[1,1,1]^T$ and shrink with μ. Here it is described in *intercept form* so it could be clearly seen in figure 6 the various variable μ means various axis intercepts of the (hyper)plane with fixed direction. This conventional optimisation problem is to find maximum possible the (hyper) plane $\boldsymbol{\Lambda}$ can extent within the feasible region.

In Flötteröd- Rohde example, the feasible region form by $\boldsymbol{\Gamma_0}, \boldsymbol{\Gamma_1}, \boldsymbol{\Gamma_2} \dots \boldsymbol{\Gamma_6}$ is symmetry (with respect to (hyper) plane $q_2$-$q_3$=0 ); therefore, there is 2 optimal solutions, which are: $\mathbf{q}^*_{FM-A}$=[0, 600, 200]$^T$ and $\mathbf{q}^*_{FM-B}$=[0, 200, 600]$^T$. In Jabari (2016), the *simplex method* is used; hence, (temporary and final) solutions only move from and to corners of the polyhedral.

(a) (b)

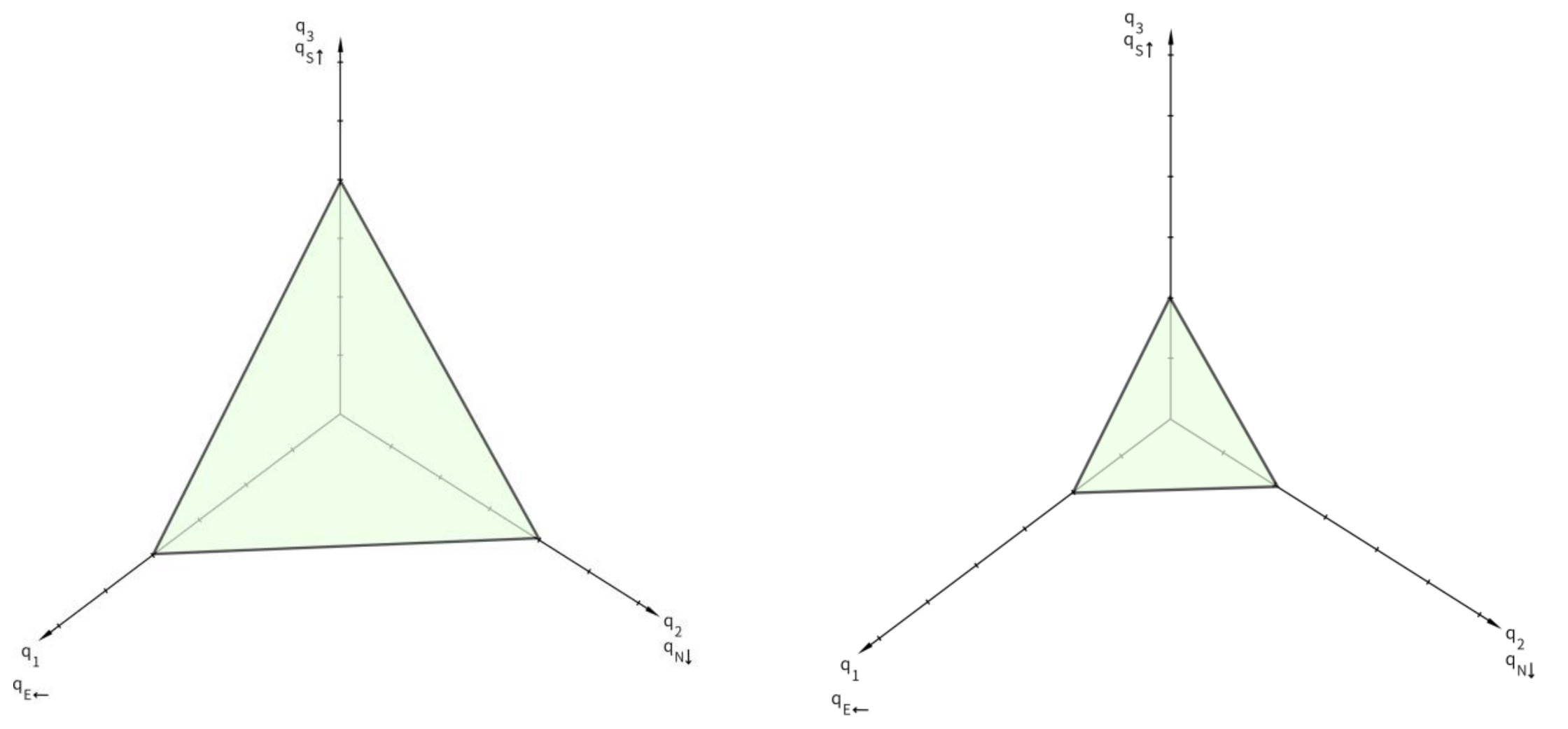

**Fig 6** Plane with respect to with different magnitude of objective function. (a) 400 (b) 200.

### 3.6 Greedy Algorithm

A greedy algorithm is proposed in Jabari (2016). Two of them has same solutions as flow maximisation: $\mathbf{q}^*_{\mathrm{Gr-A}} = \mathbf{q}^*_{\mathrm{FM-A}}$=[0, 600, 200]$^{\mathrm{T}}$ and $\mathbf{q}^*_{\mathrm{Gr-B}}$=$\mathbf{q}^*_{\mathrm{FM-B}}$=[0, 200, 600]$^{\mathrm{T}}$ Take $\mathbf{q}^*_{\mathrm{Gr-A}}$ , which comes from permutation $\mathcal{P}_{\mathrm{A}}$ in Table 1, for example, the computation steps and its geometry interpretation are:

**Step1**: $q_{\mathrm{N}} \equiv q_2$ is bounded by demand at A1. But the some of the supply is used so the residual of supply should be calculated. The residual is shown in dotted lines in Fig 7.

**Step 2a**: In this step, it is aim at finding flow from south, or $q_{\mathrm{S}} \equiv q_3$. Thus, the greedy behaviour in this step has to find another point toward $q_3$. Since the dot line is reached (in this greedy algorithm, it does not know the existence of the plane it touches in this step but rather the dotted line), we stop at A2 of Fig7.

**Step 2b**: Since $q_3$ utilised part of the supply, the residual capacity is reduced again. It consumes part of the north outgoing link and part of the west outgoing link. In other words, residual with respect to $\mathbf{\Gamma_4}$ (not shown in figure) and $\mathbf{\Gamma_5}$ (shown in dashed line in Fig 7) has changed.

**Step 3**: In this step, it is aim at finding flow from east , or $q_{\mathrm{E}} \equiv q_1$. But the residual capacity is already used up, so it ended up at A2 = $\mathbf{q}^*_{\mathrm{Gr-A}}$.

**Table 1** Solutions of greedy algorithm with different permutations

| Permutation | Solution $(i_{\mathrm{E}}, i_{\mathrm{N}}, i_{\mathrm{S}})$ |
|---|---|
| $\mathcal{P}_{\mathrm{A}} =< i_{\mathrm{N}}, i_{\mathrm{S}}, i_{\mathrm{E}} >$ | $\mathbf{q}^*_{\mathrm{Gr-A}}$=[0, 600, 200]$^{\mathrm{T}}$=$\mathbf{q}^*_{\mathrm{FM-A}}$ |
| $\mathcal{P}_{\mathrm{B}} =< i_{\mathrm{S}}, i_{\mathrm{N}}, i_{\mathrm{E}} >$ | $\mathbf{q}^*_{\mathrm{Gr-B}}$=[0, 200, 600]$^{\mathrm{T}}$=$\mathbf{q}^*_{\mathrm{FM-B}}$ |
| $\mathcal{P}_{\mathrm{C}}=< i_{\mathrm{N}}, i_{\mathrm{E}}, i_{\mathrm{S}} >$ | $\mathbf{q}^*_{\mathrm{Gr-C}}$= [100, 600, 0]$^{\mathrm{T}}$ |
| $\mathcal{P}_{\mathrm{D}}=< i_{\mathrm{S}}, i_{\mathrm{E}}, i_{\mathrm{N}} >$ | $\mathbf{q}^*_{\mathrm{Gr-D}}$= [100, 0, 600]$^{\mathrm{T}}$ |
| $\mathcal{P}_{\mathrm{E}}=< i_{\mathrm{E}}, i_{\mathrm{N}}, i_{\mathrm{S}} >$ | $\mathbf{q}^*_{\mathrm{Gr-E}}$= [100, 600, 0]$^{\mathrm{T}}$ |
| $\mathcal{P}_{\mathrm{F}}=< i_{\mathrm{E}}, i_{\mathrm{S}}, i_{\mathrm{N}} >$ | $\mathbf{q}^*_{\mathrm{Gr-F}}$= [100, 600, 0]$^{\mathrm{T}}$ |

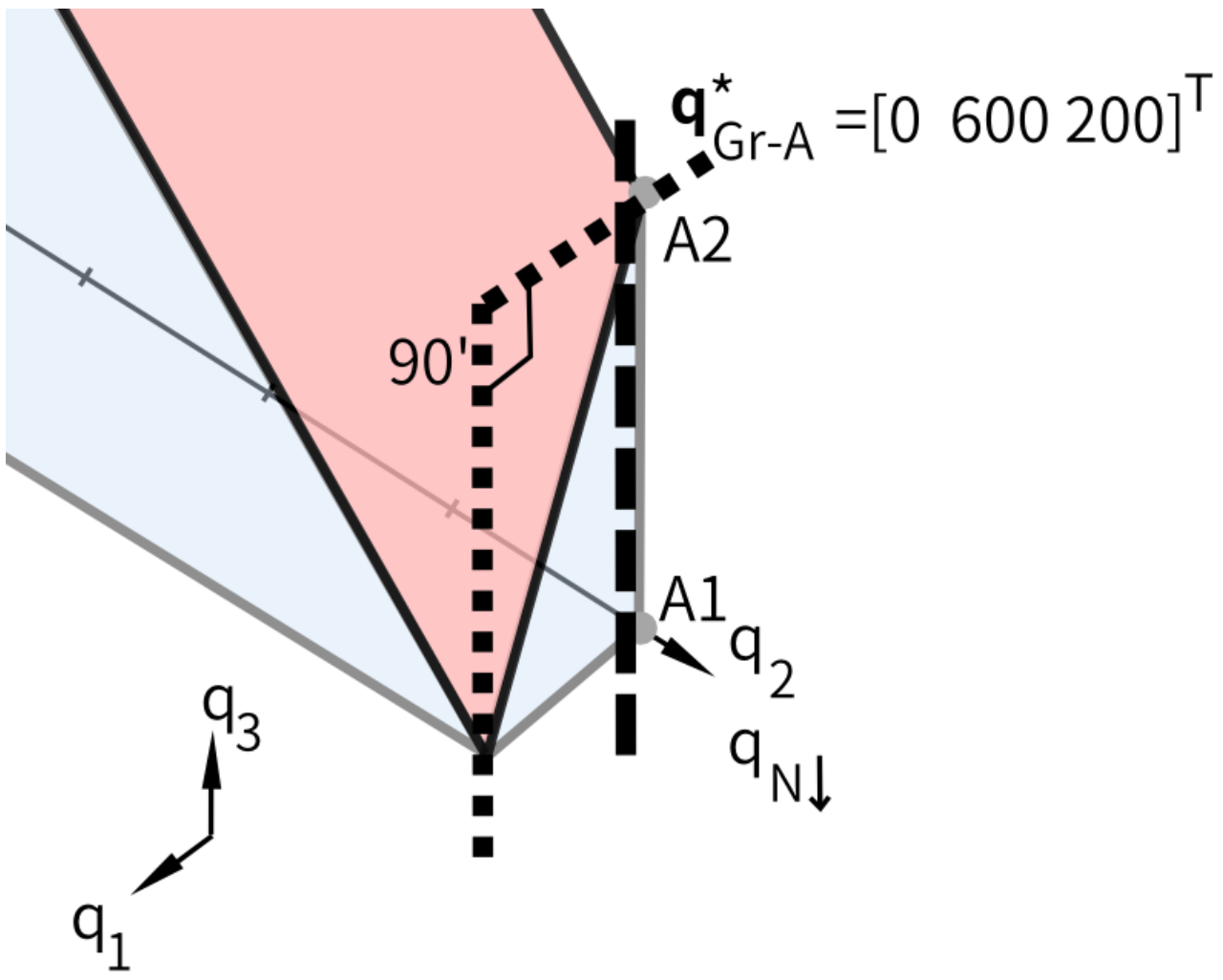

**Fig 7** The solution steps when using greedy algorithm.

### 3.7 Some Observation

From the figures presented we can observe that all the algorithms (including INM, FM-x, Gr-x) that are HFS (proofs made in Jabari(2016)) are trying to make its termination point to 'as binding' as possible. If not laying on one of the constraints, the algorithms are not terminated. This gives us some hint on relation between HFS and the behaviour of "trying to avoid slackness".

## 4. Pareto Optimal

### 4.1 Binding (Active) of the Inequality: HFS

Let us revisit the HFS necessary and sufficient condition (Jabari, 2016):

$$(\delta_i - q_i) \cdot \prod_{f(i,o)>0} (\sigma_o - q_o) = 0 \ \forall i \in I \tag{11}$$

Or it can be arranged in a less compact form:

$$(\delta_i - q_i) \cdot \prod_{f(i,o)>0} \left(\sigma_o - \sum_{i'} f(i',o) q_{i'}\right) = 0 \ \forall i \in I \tag{12}$$

It can be stated it in: For an incoming link i, its flow should be either constrained by demand or by at least one of the supply constraint which has nonzero turning ratio from i. It is equivalent to the following condition.

(Modified version of) HFS necessary and sufficient condition $\forall i \in I$:

$$(\delta_i - q_i) = 0 \ \text{ or } \ \prod_{f(i,o)>0} \left(\sigma_o - \sum_{i'} f(i',o) q_{i'}\right) = 0 \tag{13}$$

It can be added with slack variables $s_i$ $\forall i \in \mathcal{I}$ and $s_o$ $\forall o \in \mathcal{O}$ and thus arranged as:

(Another equivalent formulation) HFS necessary and sufficient condition is when objective function $\upsilon$ is found to be $\upsilon^*$=0 when solving:

$$\min \quad \upsilon = \sum_i \left(s_i \cdot \prod_{f(i,o)>0} s_o\right) \tag{14}$$

$$\begin{aligned} \text{s.t.}\ \Gamma_i': \quad & q_i + s_i = \delta_i & \forall i \in \mathcal{I} \\ \Gamma_O': \quad & \sum_{i'} f(i',o) q_{i'} + s_o = \sigma_o & \forall o \in \mathcal{O} \\ \Psi'_i: \quad & s_i \geq 0\ \forall i \in \mathcal{I} \\ \Psi'_O: \quad & s_o \geq 0\ \forall o \in \mathcal{O} \\ \Gamma_0: \quad & q_i \geq 0\ \forall i \in \mathcal{I} \end{aligned} \tag{15}$$

The problem ( eqn(14) and (15) as a whole) has feasible region as a convex set in eqn(15), but with a nonlinear objective function in eqn(14). It forms an augmented space in $\mathbb{R}^{2\cdot|\mathcal{I}|+|\mathcal{O}|}$ (instead of just in $\mathbb{R}^{|\mathcal{I}|}$)

Its interpretation is as follows: At least one of the nontrivial constraints is not slack for each incoming link. So if HFS holds the solution must be on the nonzero boundary of feasible polyhedral at the same time it cannot have a better flow performance in the neighbouring feasible region; and the better flow performance means when moving from point E to point F, then for any incoming link i with higher flow (smaller slackness $s_i$ for i) does not depreciate the performance of another incoming link i'(greater slackness $s_i{}'$ for i').

**4.2 Problem Nature: Pareto Optimal**

Let us investigate what happens if moving from points on $\upsilon^*$ =0 to other points $\bar{\upsilon}$ =0 and $\mathbf{q}^* \neq \bar{\mathbf{q}}$. It will be shown later that it is in fact Pareto optimal. But we have to define what 'better off' is: $\tilde{\mathbf{q}}$ is said to be 'better off' than original $\mathbf{q}^*$ **(**or, $\tilde{\mathbf{q}}$ is said to dominates over $\mathbf{q}^*$), denoted in $\tilde{\mathbf{q}} \succ \mathbf{q}^*$.

$$\tilde{\mathbf{q}} \succ \mathbf{q}^* \text{ if } \mathbf{q}^* \equiv [q_1^* \ldots q_i^* \ldots]^T,\ \tilde{\mathbf{q}}^* \equiv [\tilde{q}_1 \ldots \tilde{q}_i \ldots]^T,\ \tilde{q}_i > q_i^*,\ \forall i \in \mathcal{I}. \tag{16}$$

Moving to a new point without changing the objective value, then one of the incoming link must change at least one of its relevant slackness ($s_i$ and $s_o$ for all $f(i,o) > 0$) other than the zero slack variable. Then one of the following conditions (A and B) holds:

**(Condition A)** $\exists i$ $s^*{}_i = \bar{s}_i$ =0 and $\exists$ $o \in \{o\text{:}f(i,o)>0\}$ $\bar{s}_o \neq s^*{}_o$ (otherwise, if $\forall i$ $s^*{}_i = \bar{s}_i$ and $\bar{s}_o = s^*{}_o$, it implies $\mathbf{q}^* = \bar{\mathbf{q}}$ , but we want to investigate the behaviour of $\mathbf{q}^* \neq \bar{\mathbf{q}}$) Then it can be further categorised into one of the following conditions:

**(A1)** $s^*_o < \bar{s}_o$: $s^*_i = \bar{s}_i$ implies $q^*_i = \bar{q}_i$. From $f(i,o)q_i + \sum_{i''} f(i'',o)q_{i''} + s_o = \sigma_o$, we can see that $\sum_{i''} f(i'',o)q^*_{i''} > \sum_{i''} f(i'',o)\bar{q}_{i''}$. And this implies that one of the other incoming flow $q_{i'}$ makes new state $\bar{\mathbf{q}}$ no 'better off' since $\exists i'$ $q^*_{i'} > \bar{q}_{i'}$.

**(A2)** $s^*_o > \bar{s}_o$: $s^*_i = \bar{s}_i$ implies $q^*_i = \bar{q}_i$. And in $f(i,o)q_i + \sum_{i''} f(i'',o)q_{i''} + s_o = \sigma_o$ $\exists i'$ $q^*_{i'} < \bar{q}_{i'}$ It implies that this other flow $q_{i'}$ is restricted by supply constraint instead of demand constraint (otherwise its value would stay the same). We then categorise $i'$ into condition B (Rename $i'$ as i in Condition B).

**(Condition B)** $\exists i$ such that $\exists$ $o \in \{o{:}f(i,o)>0\}$ $\bar{s}_o = s^*_o = 0$. Thus, either $s^*_i \neq \bar{s}_i$ (categorised as condition B1 or B2) or $s^*_{o'} \neq \bar{s}_{o'}$ (condition B3 or B4) or both

**(B1)** $s^*_i < \bar{s}_i$. $q^*_i > \bar{q}_i$, so $\bar{q}_i$ makes new state $\bar{\mathbf{q}}$ is not 'better off'.

**(B2)** $s^*_i > \bar{s}_i$. Since $s_o$ is fixed and $q^*_i < \bar{q}_i$, from $f(i,o)q_i + \sum_{i''} f(i'',o)q_{i''} + s_o = \sigma_o$ it can be seen that $\sum_{i''} f(i'',o)q^*_{i''} > \sum_{i''} f(i'',o)\bar{q}_{i''}$ And this implies that one of the other incoming flow $q_{i'}$ become no 'better off' since $\exists i'$ $q^*_{i'} > \bar{q}_{i'}$.

**(B3)** $s^*_{o'} < \bar{s}_{o'}$ From $f(i,o')q_i + \sum_{i''} f(i'',o')q_{i''} + s_{o'} = \delta_{o'}$, we can see that either $q^*_i > \bar{q}_i$ or $\sum_{i''} f(i'',o)q^*_{i''} > \sum_{i''} f(i'',o)\bar{q}_{i''}$ or both. $q^*_i > \bar{q}_i$ implies new state is no 'better off' and $\sum_{i''} f(i'',o)q^*_{i''} > \sum_{i''} f(i'',o)\bar{q}_{i''}$ implies some of the flow in the summation, $q_{i'}$, become no 'better off' since $\exists i'$ $q^*_{i'} > \bar{q}_{i'}$.

**(B4)** $s^*_{o'} > \bar{s}_{o'}$. From $f(i,o')q_i + \sum_{i''} f(i'',o')q_{i''} + s_{o'} = \sigma_{o'}$, it can be inferred that either $f(i,o')q^*_i < f(i,o')\bar{q}_i$ or $\sum_{i''} f(i'',o')\bar{q}_{i''} < \sum_{i''} f(i'',o')q^*_{i''}$ (or both), then it can be further categorised as follows:

(i) $q^*_i < \bar{q}_i$. In the still binding constraint, $\bar{s}_o = s^*_o = 0$. From $f(i,o)q_i + \sum_{i'''} f(i''',o)q_{i'''} = \sigma_o$, we can see that for some i'''' in the summation: $q^*_{i''''} > \bar{q}_{i''''}$, so the new state is no 'better off'.

(ii) $\sum_{i''} f(i'', o')\bar{q}_{i''} < \sum_{i''} f(i'', o') q^*_{i''}$. There exist some i" in the summation having the property $\bar{q}_{i''} < q^*_{i''}$. This incoming link $i''$ must be originally bounded by supply constraint actively (it must not be originally bounded by demand constraint; otherwise it would be infeasible in new state). Or, flow of link $i''$ is involved in another binding on outgoing link $o'''$ as: $f(i'', o''')q_{i''} + \sum_{\iota} f(\iota, o''')q_{\iota} + 0 = \sigma_{o'''}$. Then, for some flow $\iota$, $\bar{q}_{\iota} > q^*_{\iota}$ which makes no 'better off' on new state.

All of the conditions show that new state $\bar{\mathbf{q}}$ which also minimises $\upsilon$ cannot be 'better off'. We thus can conclude that:

Define Y as the polyhedral form by $\Gamma_0, \Gamma_i(\forall i \in I)$ and $\Gamma_O(\forall o \in \mathcal{O})$. And define P(Y) as the subset of Y, then the HFS satisfies: $P(Y) = \{ \mathbf{q}^* \in Y \colon \{ \tilde{\mathbf{q}} \in Y : \tilde{\mathbf{q}} \succ \mathbf{q}^*, \tilde{\mathbf{q}} \neq \mathbf{q}^* \} = \emptyset \}$, which is Pareto optimal (Definition of Pareto optimal can be seen in, for instance, Kalai and Samet,1985)

In other words, the objective function $\upsilon^* = 0$ we found is always Pareto optimal. Or, $\upsilon = 0$ (another formulation of HFS) is a sufficient condition for Pareto optimal. As shown in Fig 7, all solutions from INM, flow maximisation, and greedy algorithm stands in Pareto frontier, which is shown in area with wave mark (including the boundary).

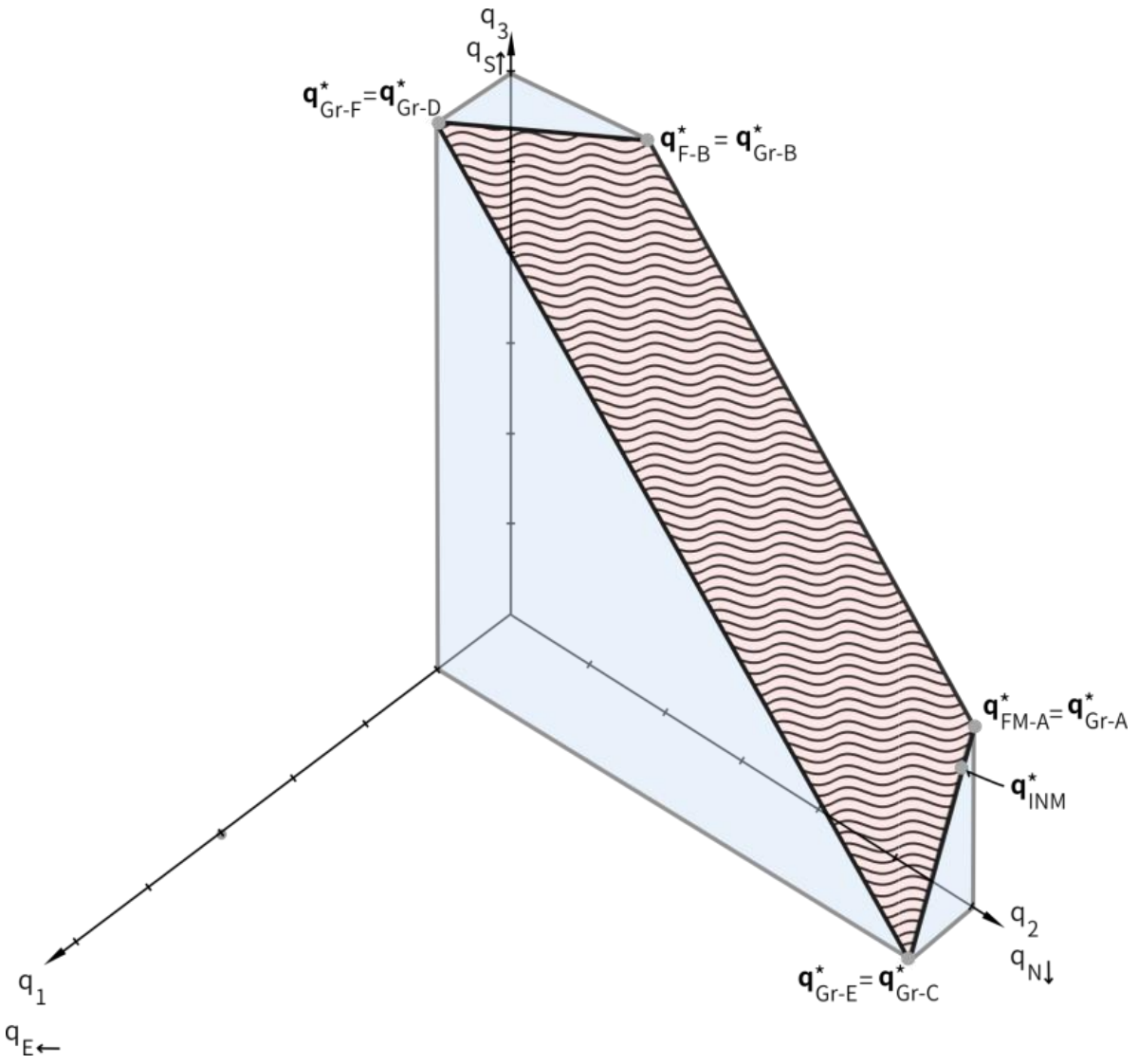

**Fig 8** All HFS lies on the Pareto frontier, the area in wave mark (including its boundary).

## 5. Discussion

### 5.1 The 'Sub-optimality' and optimality

Jabari (2016) points out that the solution from INM is suboptimal, it is true in flow maximising sense because total flow of INM (783.3 veh/hr) is less than results from FM-A, FM-B, Gr-A, and Gr-B (800 veh/hr). However, as we formulate the problem in HFS condition in an optimisation problem, INM reaches optimality (in augmented space which includes slack variables) because the optimised objective function is reached at zero ( $\upsilon^* = 0$). No matter INM, Flow Maximisation, or Greedy Algorithm, they all reaches (Pareto) optimal, but might stand on different positions of Pareto frontier.

In this perspective, flow to reach maximisation is not necessary if one is trying to find Pareto optimal only, which can be reached if HFS is reached.

### 5.2 Flexibility

Method of INM is more flexible than flow maximisation method (and the greedy algorithm). The solution of flow maximisation would only be appeared on the corners of polyhedral because solved from simplex method (Jabari, 2016). Even taking the permutation differently in greedy algorithm one can still derive corner solutions only, leading to possibly some of incoming flow zero, which might not be observed in practice.

The INM, on the other hand, is more flexible because the parameter of merging weight can be used to tune the flows not to the corners. Geometrically, merging weight $\boldsymbol{\alpha}$ provide a line in flow space to intersect with the binding (hyper) plane. Thus, the solutions will not be constrained only at corners.

## 6. Summary

For 3 or less incoming links, the constraints can be depicted in a 3(or less)-dimensional diagram. The nature of Flötteröd- Rohde example (Flötteröd and Rohde, 2011) happened to fall into this category and allows us to see some insights from geometric viewpoint.

Although in the beginning the analysis of drawing diagram by hand it is restricted to up to 3 dimensional, analysing slack variables thereafter does not have the restriction to low dimension. In higher dimensions in general (which means more than 3 incoming links), the holding-free solution lies on frontiers of the polyhedral, so finding the HFS is rather a problem of Pareto frontier-finding instead of ordinary optimisation.

In this paper, only simple supply/demand constraint is investigated. The next phase of the research would be including the *internal supply constraint* and the relevant INM with *internal supply constraint* (INMC) (Flötteröd and Rohde, (2011)) into geometry interpretation.

**Reference**

1. Corthout, R., Flötteröd, G., Viti, F., Tampère, C.M., 2012. Non-unique flows in macroscopic first-order intersection models. Transportation Research Part B: Methodological 46 (3), 343–359.

2. Flötteröd, G., Rohde, J., 2011. Operational macroscopic modeling of complex urban road intersections. Transportation Research Part B 45, 903-922.

3. Gibb, J., 2011. Model of traffic flow capacity constraint through nodes for dynamic network loading with queue spillback. Transportation Research Record: Journal of the Transportation Research Board 2263 (-1), 113–122.

4. Jabari, S.E., 2016. Node modeling for congested urban road networks. Transportation Research Part B 91, 229–249.

5. Kalai, E., and Samet D., 1985. Unanimity games and Pareto optimality. International Journal of Game Theory 14, 41-50.

6. Lebacque, J., Khoshyaran, M., 2005. First order macroscopic traffic flow models: intersections modeling, network modeling. In: Proceedings of the 16th ISTTT.

7. Smits, E.S., Bliemer, M. C.J, Pel, A. J., van Arem, B., 2015. A family of macroscopic node models. Transportation Research Part B: Methodological 74 (1) 20–39

8. Tampère, C.M., Corthout, R., Cattrysse, D., Immers, L.H., 2011. A generic class of first order node models for dynamic macroscopic simulation of traffic flows. Transportation Research Part B: Methodological 45 (1), 289–309.

9. Transportation Research Board, National Research Council, 2010. Highway Capacity Manual 2010. Washington DC.